\documentclass[12pt]{amsart}

\usepackage{geometry}                

\usepackage[all]{xy}

\usepackage{graphicx}
\usepackage{amssymb}
\usepackage{epstopdf}

\usepackage{amsmath, amsthm}
\usepackage{mathrsfs}
\usepackage{url}

\usepackage{color}
\definecolor{MyRed}{rgb}{0.6,0,0}

\newcommand{\tg}{\pmb{\top}}
\newcommand{\rr}{\mathbb{R}}
\newcommand{\zz}{\mathbb{Z}}
\newcommand{\qq}{\mathbb{Q}}
\newcommand{\cc}{\mathbb{C}}

\newcommand{\ccc}[1]{\mathbb{C}^{#1}}

\newcommand{\zplus}{\mathbb{Z}_{>0}}

\newcommand{\rrnn}{\mathbb{R}_{\geq 0}}
\newcommand{\zznn}{\mathbb{Z}_{\geq 0}}
\newcommand{\comproj}[1]{\mathbb{CP}^{#1}}
\newcommand{\reproj}[1]{\mathbb{RP}^{#1}}

\newcommand{\lie}{\mathfrak{g}}

\newcommand{\liek}{\mathfrak{k}}
\newcommand{\lieb}{\mathfrak{b}}
\newcommand{\lieu}{\mathfrak{u}}
\newcommand{\liet}{\mathfrak{t}}
\newcommand{\lieq}{\mathfrak{q}}

\newcommand{\chamber}{\mathfrak{t}_{+}^{*}}

\newcommand{\lied}{\mathcal{L}}
\newcommand{\iprod}[2]{\langle #1, #2 \rangle}

\newcommand{\Hom}{\operatorname{Hom}}

\newcommand{\pr}{\operatorname{pr}}
\newcommand{\real}{\operatorname{Re}}
\newcommand{\imaginary}{\operatorname{Im}}
\newcommand{\ii}{\sqrt{-1}}
\newcommand{\func}[3]{#1 \colon #2 \to #3}

\newcommand{\deff}[1]{\emph{\bf #1}}

\theoremstyle{plain}
\newtheorem{theorem}{Theorem}[section]

\newtheorem{corollary}[theorem]{Corollary}
\newtheorem{lemma}[theorem]{Lemma}
\newtheorem{proposition}[theorem]{Proposition}

\theoremstyle{definition}

\newtheorem{remark}[theorem]{Remark}

\newtheorem{definition}[theorem]{Definition}

\subjclass[2000]{53D20 (14L24 53C55)}

\title{A convexity theorem for the real part of a Borel invariant subvariety}
\author{Timothy E. Goldberg}
\address{Department of Mathematics \\ Cornell University \\ Ithaca, NY 14850-4201}
\email{goldberg@math.cornell.edu}
\urladdr{http://www.math.cornell.edu/~goldberg}
\date{September 9, 2008}                                           
\thanks{The author was partially supported by National Science Foundation Grant DMS--0300172.} 

\begin{document}

\maketitle

\begin{abstract}
M. Brion proved a convexity result for the moment map image of an irreducible subvariety of a compact integral K\"ahler manifold preserved by the complexification of the Hamiltonian group action. V. Guillemin and R. Sjamaar generalized this result to irreducible subvarieties preserved only by a Borel subgroup. In another direction, L. O'Shea and R. Sjamaar proved a convexity result for the moment map image of the submanifold fixed by an antisymplectic involution. Analogous to Guillemin and Sjamaar's generalization of Brion's theorem, in this paper we generalize O'Shea and Sjamaar's result, proving a convexity theorem for the moment map image of the involution fixed set of an irreducible subvariety preserved by a Borel subgroup.
\end{abstract}

\tableofcontents

\section*{Introduction}

At the end of their 1982 paper, \cite{gs}, Guillemin and Sternberg gave a description of the moment map image of an integrable symplectic manifold with a Hamiltonian action of a compact group in terms of certain of the highest weights of a maximal torus of the group. In 1987, \cite{brion}, Brion expanded this technique and applied it to certain algebraic subvarieties of the manifold, proving a convexity theorem for their moment images, and also describing each moment polytope in terms of the \textbf{highest weight polytope}, defined in Section~\ref{sec2} below. These methods have proved very useful. They were central to Guillemin and Sjamaar's 2006 generalization of Brion's theorem, \cite{guill}, and in proving the projective case of O'Shea and Sjamaar's theorem from 2000, \cite[Section 6]{oshea}. Not surprisingly, the highest weight polytope is the main tool in this paper as well. It is well--known that describing the moment polytope is often at least as difficult as proving that the moment image is a convex polytope in the first place. We are fortunate to be able to make some descriptions here.

In Section~\ref{sec1} below, we lay out some of the technical context of the results of Brion, Guillemin--Sjamaar, and O'Shea--Sjamaar mentioned above. We then proceed to describe their results in more detail, leading to a statement of the main theorem of this article. In Section~\ref{sec2}, we describe the main tool in analyzing the moment map image in this context, the highest weight polytope. The proof of the main theorem comes in Section~\ref{sec3}. Section~\ref{sec4} contains some easy but interesting corollaries to the main theorem, and Section~\ref{sec5} describes a specific example in which the main theorem can be applied.

\section{Background} \label{sec1}

Let $M$ be an integral K\"ahler manifold which is compact and connected, and let $L$ be a holomorphic line bundle over $M$ whose Chern class is the cohomology class represented by $M$'s K\"ahler form, $\omega$. Then there is a Hermitian metric on $L$ with metric connection $\nabla$ whose curvature form satisfies $\operatorname{curv} \nabla = \frac{1}{2 \pi \ii} \, \omega$. (See \cite{wells}, for instance.) Let $G$ be a compact and connected Lie group that acts in a K\"ahlerian fashion on $(M,L)$, so that $G$ acts by complex linear bundle automorphisms. The group of all line bundle automorphisms preserving the holomorphic structure of $L$ is a complex Lie group, so the action of $G$ lifts to a holomorphic action of the complexification $G^{\cc}$ of $G$ on $(M,L)$. The Kodaira embedding theorem implies that $M$ can be embedded in some complex projective space as a closed and complex algebraic variety, on which the action of $G^{\cc}$ is algebraic.

\begin{remark} Note that for any $p \in M$, the $G^{\cc}$-orbit through $p$ is the image of the algebraic map $G^{\cc} \to M$, $g \mapsto g \cdot p$, and so by Chevalley's Theorem is a constructible set. This means that its Zariski closure and its closure in the topology of the manifold coincide. (See Corollary 2 in Section I.8 and Corollary 1 in Section I.10 of \cite{mumford}.) The same is true for an orbit of any algebraic subgroup of $G^{\cc}$, such as a Borel subgroup. Finally, note all complex subvarieties of $M$ that are closed in the manifold topology are also constructible, and hence closed in the Zariski topology as well. Therefore, for all of these sets, there is no distinction between ``closed'' and ``Zariski--closed''.
\end{remark}

Because $G$ preserves the Hermitian structure of $L$, it also preserves the symplectic form $\omega$ of $M$ given by its K\"ahler structure, so $G$ acts by symplectomorphisms. Furthermore, the action of $G$ on $M$ is Hamiltonian, with moment map $\func{\Phi}{M}{\lie^*}$ obtained as follows. Let $s$ be any global smooth section of $(M,L)$. Each $\xi \in \lie$ acts on $s$ in two ways: Lie differentiation $\lied(\xi)$, and covariant differentiation $\nabla \left(\xi_M\right)$ in the direction of the fundamental vector field $\xi_M$ on $M$. In \cite[Theorem 4.3.1]{kostant}, Kostant showed that their difference $\lied(\xi)-\nabla(\xi_M)$ is multiplication by an imaginary-valued function on $M$. Hence we can define a real linear map $\lie \to (M \to \rr)$, $\xi \mapsto \phi^{\xi}$ by \[\phi^{\xi} = \frac{1}{2 \pi \ii} \left( \lied(\xi) - \nabla(\xi_M) \right)\] for each $\xi \in \lie$. Then $\Phi$ is defined by the equation $\Phi(x) \, \xi = \phi^{\xi}(x)$, for all $\xi \in \lie$ and $x \in M$. It can be shown that $\Phi$ satisfies the properties of a moment map. This description comes from \cite[Section 2]{guill}. The assumption that the moment map here is not arbitrary, but is intimately connected to the actions of $G$ on both $M$ and $L$, is extremely fruitful.

Suppose we have involutions $\func{\gamma}{G}{G}$, $\func{\tau}{M}{M}$, and $\func{\beta}{L}{L}$ such that $\gamma$ is a smooth group homomorphism, $\tau$ is antiholomorphic and antisymplectic, and $(\tau, \beta)$ is an involutive (real) bundle automorphism on $(M,L)$ which is complex antilinear on fibers and which preserves the covariant derivative $\nabla$ on $L$. Then $\func{\gamma}{G}{G}$ induces linear involutions $\func{\gamma}{\lie}{\lie}$ and $\func{\gamma}{\lie^*}{\lie^*}$, defined in the obvious way. We further require two compatibility conditions regarding the involutions and the Hamiltonian action of $G$ on $M$. We assume the properties of
\begin{enumerate}
\item[(1)] distribution, i.e. $\tau (g \cdot m) = \gamma(g) \cdot \tau(m)$ for all $g \in G, m \in M$; and
\item[(2)] anti-equivariance, i.e. $\Phi\left( \tau(m) \right) = - \gamma \left( \Phi(m) \right)$.
\end{enumerate}

There are two obvious ways to extend $\func{\gamma}{\lie}{\lie}$ to an involution on its complexification $\lie^{\cc} := \lie \otimes_{\rr} \cc$ --- holomorphically or antiholomorphically. The antiholomorphic is more useful for our purposes. Define $\func{\sigma}{\lie^{\cc}}{\lie^{\cc}}$ by $\sigma(\xi) = \gamma(\real \xi) - \ii \, \gamma(\imaginary \xi)$ for all $\xi \in \lie^{\cc}$, where $\real$ and $\imaginary$ denote the real and imaginary parts with respect to the decomposition $\lie^{\cc} = \lie \oplus \ii \lie$. The antiholomorphic Lie algebra involution $\sigma$ lifts to an antiholomorphic Lie group involution on $G^{\cc}$, which we will also denote by $\sigma$. In \cite[Proposition 5.5]{oshea}, it is proved that under the compatibility conditions described above, the fixed set of $G^{\cc}$ under this antiholomorphic involution preserves the fixed set of $M$ under $\tau$, and this is the key property we need.

Recall that a linear involution on a vector space is diagonalizable, and has eigenvalues both or one of $\pm 1$. Let $\liek$ and $\lieq$ denote the $1$ and $-1$-eigenspaces of $\func{\gamma}{\lie}{\lie}$, respectively. We can identify $\liek^*$ and $\lieq^*$ with the annihilators of $\lieq$ and $\liek$, respectively, and obtain a decomposition $\lie^* = \liek^* \oplus \lieq^*$ which is also the decomposition of $\lie^*$ into eigenspaces of $\func{\gamma}{\lie^*}{\lie^*}$. From the definition of $\sigma$, we see that the $1$-eigenspace of $\lie^{\cc}$ under this involution is exactly $\liek \oplus \ii \lieq$, and the $(-1)$-eigenspace is $\lieq \oplus \ii \liek$.

As usual, we denote fixed sets of the actions of these involutions and these groups by superscripts. One well--known reason for requiring $\tau$ to be antisymplectic is that the submanifold $M^{\tau}$ of $M$ is then Lagrangian.

\begin{proposition}
\label{lagprop}
Let $(M, \omega)$ be a symplectic manifold and $\tau$ be an antisymplectic involution on $M$. Then $M^{\tau}$ is a Lagrangian submanifold of $M$.
\end{proposition}

\begin{proof}
That $M^{\tau}$ is a submanifold of $M$ is well--known, in particular since $\tau$ defines an action of the compact group $\zz / 2 \zz$ on $M$, so that $M^{\tau} = M^{\zz/2 \zz}$. To prove that $M^{\tau}$ is Lagrangian, it suffices to show that each tangent space of $M^{\tau}$ is a Lagrangian subspace of the corresponding tangent space of $M$. Let $x \in M$, let $\Omega := \omega_x$, and let $\tilde{\tau} := \tg_x \tau$ be the derivative of $\tau$ at $x$, which is a linear involution on the tangent space $\tg_x M$. Let $\tg_x M = V^+ \oplus V^-$ be the decomposition of $V$ into $\pm 1$ eigenspaces with respect to $\tilde{\tau}$. Set $V:= \tg_x M^{\tau}$, and note that $V = \left( \tg_x M \right)^{\tilde{\tau}} = V^+$. Since $\tau^* \omega = -\omega$, we have $\Omega(u,v) = -\Omega \left( \tilde{\tau}(u), \tilde{\tau}(v) \right)$ for all $u,v \in V$. As usual, we denote symplectic complements by a superscript $\Omega$.

Notice that $\Omega$ vanishes on $V^{\pm}$, for if $v,v \in V^{\pm}$ then $\Omega(u,v) = -\Omega \left( \tilde{\tau}(u), \tilde{\tau}(v) \right) = - \Omega(\pm u, \pm v) = -\Omega(u,v)$, so $\Omega(u,v)=0$. Since $V  =V^+$, it follows that $V \subset V^{\Omega}$, so $V$ is isotropic.

Let $u,v \in V^{\Omega}$, and let $u=u^+ + u^-$ and $v=v^+ + v^-$ denote their decompositions according to $\tg_x M = V^+ \oplus V^-$. Again since $V=V^+$, we know $\Omega(u,v) = \Omega(u, v^+) + \Omega(u,v^-) = \Omega(u,v^-) = \Omega(u^+,v^-) + \Omega(u^-,v^-) = \Omega(u^+,v^-)$. By the antisymmetry of $\Omega$, we also have $\Omega(u,v) = \Omega(u^-,v^+)$. Therefore $\Omega(u,v) = \Omega(u^+ + u^-,v^+ + v^-) = \Omega(u^+,v^+) + \Omega(u^+,v^-) + \Omega(u^-,v^+) + \Omega(u^-,v^-) = \Omega(u^+,v^-) + \Omega(u^-,v^+) = 2 \, \Omega(u,v)$, so $\Omega(u,v) = 0$. Hence $V^{\Omega} \subset \left( V^{\Omega} \right)^{\Omega} = V$, so $V$ is coisotropic. Thus $V$ is Lagrangian.
\end{proof}

Let $T$ be a maximal torus of $G$ with Lie algebra $\liet$, and suppose it is preserved by $\gamma$. (As described in Appendix B of \cite{oshea}, such a torus can always be obtained by starting from a maximal torus of the submanifold $Q = \{ g  \, \gamma(g)^{-1} \mid g \in G \}$ of symmetric elements of $G$.) Choose a closed positive Weyl chamber $\chamber \subset \liet^*$. Embed $\liet^*$ as a vector subspace of $\lie^*$ in the usual way, using the real version of the root space decomposition of $\lie^{\cc}$. For any subset $A \subset M$, we let $\Delta(A) := \Phi(A) \cap \chamber$. Notice that if $m \in M^{\tau}$, then $\gamma\left(\Phi(m)\right) = -\Phi(m)$, so $\Phi(m) \in \lieq^*$. Thus $\Phi(M^{\tau}) \subset \lieq^*$. The main result of \cite{oshea} was the following essential converse. The proof required that the torus $T$ and the positive Weyl chamber $\chamber$ be chosen so as to be ``compatible'' with the involutions in a certain sense, as detailed in \cite[Section 3]{oshea}.

\begin{theorem} \label{osheathm2}
Suppose $T$ and $\chamber$ are ``compatible'' with the involutions. Then $\Delta(M^{\tau}) = \Delta(M) \cap \lieq^*$.
\end{theorem}

Later, to the current author, Sjamaar suggested and outlined the following corollary and proof. It generalizes the result of Theorem~\ref{osheathm2}, doing away with the full compatibility requirements on $T$ and $\chamber$.

\begin{corollary}[due to Sjamaar] \label{coro1} The equation $\Phi(M^{\tau}) = \Phi(M) \cap \lieq^*$ holds. Therefore, Theorem~\ref{osheathm2} is true for any choice of $T$ and $\chamber$ such that $T$ is $\gamma$--invariant.
\end{corollary}

\begin{proof}
In Example 2.9 of \cite{oshea}, the authors describe how if $\lambda \in \lieq^*$, a compatible involution $\alpha$ on the Hamiltonian $G$-manifold $G \cdot \lambda$, the coadjoint orbit through $\lambda$, is given by $\alpha := -\gamma$. In Proposition 2.3 of \cite{oshea}, they prove that $\left( G \cdot \lambda \right)^{\alpha} = G \cdot \lambda \cap \lieq^* = G^{\gamma} \cdot \lambda$.

Now let $\lambda \in \Phi(M) \cap \lieq^*$. Since $\chamber$ is a fundamental domain for the action of $G$ on $\lie^*$, there is some $g \in G$ such that $g \cdot \lambda \in \chamber$. Put $\lambda' = g \cdot \lambda$, and note that $\lambda \in G \cdot \lambda' \cap \lieq^*$. By the previous paragraph, there is some $k \in G^{\gamma}$ such that $\lambda = k \cdot \lambda'$, so $\lambda' = k^{-1} \cdot \lambda$. Since $k^{-1} \in G^{\gamma}$, we have $\gamma(\lambda') = \gamma \left( k^{-1} \cdot \lambda \right) = \gamma(k^{-1}) \cdot \gamma(\lambda) = k^{-1} \cdot (-\lambda) = - (k^{-1} \cdot \lambda) = - \lambda'$, so $\lambda' \in \lieq^*$. Because $\Phi$ is $G$-equivariant, if $\lambda = \Phi(x)$, then $\lambda' = k^{-1} \cdot \lambda = k^{-1} \cdot \Phi(x) = \Phi(k^{-1} \cdot x)$, so $\lambda' \in \Phi(M)$. Therefore $\lambda' \in \Phi(M) \cap \lieq^* \cap \chamber = \Phi(M^{\tau}) \cap \chamber$. So there is some $y \in M^{\tau}$ with $\Phi(y) = \lambda'$, which means \[ \lambda = k \cdot \lambda' = k \cdot \Phi(y) = \Phi(k \cdot y) .\] Because $k \in G^{\gamma}$ and $y \in M^{\tau}$, we have $\tau(k \cdot y) = \gamma(k) \cdot \tau(y) = k \cdot y$, so $k \cdot y \in M^{\tau}$ and $\lambda \in \Phi(M^{\tau})$. Thus $\Phi(M) \cap \lieq^* \subset \Phi(M^{\tau})$. The other inclusion was shown above.
\end{proof}

Theorem~\ref{osheathm2} and Corollary~\ref{coro1} and their proofs do not require the presence of the line bundle or the complex structures whose existence we have assumed. By Kirwan's convexity theorem (\cite{kir}), the set $\Delta(M)$ is a convex polytope in $\liet^*$, so $\Delta(M^{\tau})$ is the intersection of a convex polytope with a linear subspace, which means it too is a convex polytope.

In the full K\"ahler and line bundle circumstances we have defined here, O'Shea and Sjamaar also proved the following statement, \cite[Theorem 5.10]{oshea}.

\begin{theorem} \label{osheathm}
Let $X$ be a closed, irreducible, complex subvariety of $M$ preserved by $G^{\cc}$ and $\tau$, and let $Y \subset M^{\tau}$ be the closure of any nonempty component of $X_{\text{reg}} \cap M^{\tau}$, where $X_{\text{reg}}$ denotes the set of regular points in $X$. Then $\Delta(Y) = \Delta(X) \cap \lieq^*$.
\end{theorem}

The main result of Brion in \cite{brion} implies that $\Delta(X)$ is a convex polytope in $\liet^*$, so as before $\Phi(Y)$ is a convex polytope as well.

In \cite{guill}, Sjamaar and Guillemin strengthened Brion's convexity result. Let $B \subset G^{\cc}$ be the Borel subgroup determined by our choice $\chamber$ of positive Weyl chamber. 

\begin{theorem} \label{guillthm}
Let $X$ be a $B$-invariant irreducible closed subvariety of $M$. Then $\Delta(X)$ is a rational convex polytope in $\liet^*$, the closure of the set $\mathcal{C}(X)$, (defined in Section~\ref{sec2}), which is a convex polytope in the space of rational points in $\liet^*$.
\end{theorem}

Here, rational means rational with respect to the weight lattice of $T$, embedded in a particular way in $\liet^*$, which we specify later. This theorem and its proof do not involve any involutions, of course.

Our main result is a combination of Theorems~\ref{osheathm} and \ref{guillthm}.

\begin{theorem}[Main Theorem] \label{mainthm}
Suppose the Borel subgroup $B$ is preserved by the involution $\sigma$ on $G^{\cc}$. Let $X$ be a closed, irreducible, complex subvariety of $M$ preserved by both $B$ and $\tau$, and let $Y$ be the closure of any nonempty component of $X_{\text{reg}} \cap M^{\tau}$. Then \[\Delta(Y) = \Delta(X) \cap \lieq^*\] and $\Delta(Y)$ is a rational convex polytope in $\liet^*$, the closure of the set $\mathcal{C}_{\gamma}(Y)$, (defined in Section~\ref{sec2}), which is a convex polytope in the space of rational points in $\liet^*$.
\end{theorem}

Theorem~\ref{mainthm} immediately implies the following.

\begin{corollary} \label{maincor}
Suppose $B$ and $X$ are as in Theorem~\ref{mainthm}. If $X^{\tau} \cap X_{\text{reg}} \neq \emptyset$, then $\Delta(X^{\tau}) = \Delta(X) \cap \lieq^*$, and so $\Delta(X^{\tau})$ is a rational convex polytope in $\liet^*$.
\end{corollary}

Notice that all of these results are specific instances of the main idea that the real part of the moment polytope is the moment polytope of the real part.

\begin{remark}
Given an anti-holomorphic involution on $G^{\cc}$, the question of whether or not there exists an invariant Borel subgroup has been studied, for instance in \cite[Section 5]{adams}. In that paper, an involution for which the answer is yes is called \deff{principal}.
\end{remark}

\section{The highest weight polytope}
\label{sec2}

We follow Brion's approach from \cite{brion}, as was done in \cite[Section 5]{oshea} and \cite{guill}, and consider certain subsets of global holomorphic sections of $(M,L)$ and its tensor powers. We will decompose these spaces into weight spaces under the action of $T$.

Let $\Lambda=\Hom \left(T,\mathbf{U}(1) \right)$ be the weight lattice of $T$. We identify $\Lambda$ with a certain lattice, also denoted $\Lambda$, in $\liet^*$ via the map $\lambda \in \Lambda \mapsto \frac{1}{2\pi \ii} \, \tg_1 \lambda \in \liet^*$. Here $\tg_1 \lambda$ denotes the derivative of the map $\lambda$ at the identity. Put $\Lambda_+ = \Lambda \cap \chamber$, the space of dominant weights. We call a point in $\liet^*$ \deff{rational} if it is contained in a rational multiple of the weight lattice. Hence the set of rational points is $\Lambda \otimes_{\zz} \qq$. 

\begin{remark} \label{rem1}
The fact that $\gamma$ preserves $T$ implies that $\gamma$ preserves the lattice $\Lambda$. So with respect to a basis for $\liet^*$ consisting of lattice elements, $\gamma|_{\liet^*}$ can be represented by a matrix with rational entries. Since the only eigenvalues of $\gamma$ are the integers $1$ and $-1$, we conclude that there exist bases for the eigenspaces $\liet^* \cap \liek^*$ and $\liet^* \cap \lieq^*$ of $\gamma|_{\liet^*}$ consisting of {\bf rational} linear combinations of lattice elements. Therefore, for each of these eigenspaces, the rational elements in the eigenspace form a dense subset of it.
\end{remark}

Let $\Gamma(M,L)$ be the space of global holomorphic sections of $(M,L)$, and for each $r \in \zznn = \{0,1,2, \ldots\}$ let $\Gamma(M,L^r)$ be the space of global holomorphic sections of the $r$-fold tensor product of $(M,L)$ over $\cc$, \[(M,L^r) = ( M, \underbrace{L \otimes L \otimes \ldots \otimes L}_{r \text{ times}} ).\] (We consider $\Gamma(M,L^0)$ to be the space of holomorphic complex-valued functions on $M$. Since $M$ is compact, we know $\Gamma(M,L^0) \cong \cc$.) Since $T$ acts on $(M,L)$ and by extension on each $(M,L^r)$ by complex bundle automorphisms, $T$ acts on the spaces of holomorphic global sections of these bundles and in particular on the smooth sections: for any such section $s$ and any $t \in T$, the action of $t$ on $s$ is defined by $(t \cdot s)(x) := t \cdot \left[ s \left( t^{-1} \cdot x \right) \right]$ for $x \in M$.

Each $\Gamma(M,L^r)$ decomposes under the action of $T$ into weight spaces: $\Gamma(M,L^r) = \bigoplus_{\lambda \in \Lambda} \Gamma(M,L^r)_{\lambda}$. Let $S=\bigoplus_{r \in \zznn} \Gamma(M,L^r)$, and for each $r \in \zznn$ put $S_{\lambda,r} = \Gamma(M,L^r)_{\lambda}$. Then $S$ has a grading by $\Lambda \times \zznn$, $S = \bigoplus_{(\lambda,r) \in \Lambda \times \zznn}  S_{\lambda,r}$. Let $N = [B,B]$ be the unipotent radical of $B$. Then this grading of $S$ descends to a grading of the $N$-invariant elements of $S$, $S^N = \bigoplus_{(\lambda,r) \in \Lambda_+ \times \zznn} S_{\lambda,r}^N$. (Recall that any weight that appears in the weight decomposition of $S^N$ must be dominant.)

For any $B$-invariant irreducible closed complex subvariety $X$ of $M$, let $I(X)$ be the homogeneous ideal of $S$ consisting of sections that vanish identically on $X$, let $I(X)^N$ denote the set of $N$-invariant sections that vanish identically on $S$, and let $A(X)$ be the quotient $A(X) = S^N / I(X)^N$.

\begin{definition} \label{def1}
The \deff{highest weight polytope} of $X$ is the subset $\mathcal{C}(X)$ of $\Lambda \otimes \qq$ defined by $\mathcal{C}(X) := \left\{ \lambda \in \Lambda \otimes \qq \mid \text{ there exists } r \in \zplus \text{ such that } r \lambda \in \Lambda_+ \text{ and } A(X)_{r \lambda,r} \neq 0 \right\}$.
\end{definition}

As detailed in \cite{brion}, $\mathcal{C}(X)$ is indeed a convex polytope in the $\qq$-vector space $\Lambda \otimes \qq$. The specific main result of that paper is that, if $X$ is preserved by all of $G^{\cc}$, then $\Delta(X) \cap (\Lambda \otimes \qq) = \mathcal{C}(X)$ and $\Delta(X)$ is the closure of $\mathcal{C}(X)$ in $\liet^*$, so $\Delta(X)$ is a rational convex polytope. The main result of \cite{guill} is exactly that the same statements hold even if $X$ is only preserved by $B$.

To put Definition~\ref{def1} another way, an element $\lambda \in \Lambda \otimes \qq$ is contained in $\mathcal{C}(X)$ if and only if there exists $r \in \zplus$ such that $r \lambda \in \Lambda_+$ and there is a section $s \in S_{r \lambda, r}^N$ which does not vanish identically on $X$. An equivalent condition is that there exists $r \in \zplus$ such that $r \lambda \in \Lambda_+$ and the irreducible representation of $G$ with highest weight $r \lambda$ is a submodule of the $G$-module $\Gamma(M,L^r)$, and there is an element of this submodule which does not vanish identically on $X$. Accordingly, for any subset $Z$ of $M$, we make the following definition.

\begin{definition} \label{def2}
The \deff{$\gamma$-highest weight set} of $Z$ is the subset $\mathcal{C}_{\gamma}(Z)$ of $(\Lambda \otimes \qq) \cap \lieq^*$ consisting of elements $\lambda \in \lieq^*$ for which there exists $r \in \zplus$ such that the irreducible representation of $G$ with highest weight $r \lambda$ is a submodule of the $G$-module $\Gamma(M,L^r)$, and there is an element of this submodule which does not vanish identically on $Z$.
\end{definition}

\section{Proof of the main theorem}
\label{sec3}

Suppose the Borel subgroup $B$ is preserved by $\sigma$. Let $X$ be a closed, irreducible, complex subvariety of $M$ preserved by $B$ and $\tau$, and let $Y$ be the closure of any nonempty component of $X_{\text{reg}} \cap M^{\tau}$.

\begin{proposition} \label{prop1}
The equality $\mathcal{C}_{\gamma}(Y) = \mathcal{C}(X) \cap \lieq^*$ holds.
\end{proposition}

\begin{proof}
From the definition of $\mathcal{C}_{\gamma}(Y)$,  the inclusion $\mathcal{C}_{\gamma}(Y) \subset \mathcal{C}(X) \cap \lieq^*$ is immediate.

For the other direction, suppose $\lambda \in \mathcal{C}(X) \cap \lieq^*$. Then there is some $r \in \zplus$ such that $r \lambda \in \Lambda_+$ and a section $s \in S_{r \lambda,r}^N$ which does not vanish identically on $X$. Similarly to the proof of Proposition~\ref{lagprop} above, observe that $Y$ contains a Lagrangian submanifold of $X_{\text{reg}}$. By the compatibility of the complex and symplectic structures of $M$, this Lagrangian submanifold is a totally real submanifold, which implies that $Y$ is Zariski--dense in $X$. Hence any holomorphic section that vanishes on all of $Y$ must vanish on all of $X$, so $s$ cannot vanish identically on $Y$. Since $\lambda \in \lieq^*$, this means that $\lambda \in \mathcal{C}_{\gamma}(Y)$.
\end{proof}

Consider the identity component of the fixed set $(G^{\cc})^{\sigma}$ of $G^{\cc}$ under the involution $\sigma$. Note that its Lie subalgebra is $\liek \oplus \ii \lieq$. Proposition 5.5 of \cite{oshea} states that $\tau$ is equivariant under the action of this subgroup on $M$, which implies that this subgroup preserves the fixed point set $M^{\tau}$. Let $H$ denote the identity component of the ``real Borel subgroup'', $B^{\sigma}$. This group has the virtue of preserving both $X$ and $M^{\tau}$, which means it also preserves $Y$. Its Lie algebra is $\lieb^{\sigma} = (\liek \oplus \ii \lieq) \cap \lieb$.

\begin{lemma} \label{lem1}
The $\gamma$-highest weight polytope $\mathcal{C}_{\gamma}(Y)$ is the set of rational points in $\Delta(Y)$.
\end{lemma}

\begin{proof}
Let $\lambda \in \mathcal{C}_{\gamma}(Y)$. By Proposition~\ref{prop1}, this means $\lambda \in \mathcal{C}(X) \cap \lieq^*$. Then there is $r \in \zplus$ and $s \in \Gamma(M,L^r)$ such that $r \lambda \in \Lambda_+$, $s \in S_{r \lambda, r}^N$, and $s$ does not vanish identically on $X$. Because $Y$ is a closed subset of the compact space $M$, it is itself compact, so there is an element $y \in Y$ where the smooth function $\|s\|^2$ takes its maximum value on $Y$. Recall that since $s$ does not vanish on $X$, it does not vanish on $Y$, so $\|s(y)\|^2 > 0$. Because $H$ preserves $Y$, $Y$ is a union of $H$-orbits. Since the fundamental vector fields induced by elements of its Lie algebra $\lieb^{\sigma}$ are tangent to the $H$-orbit through the point at which the vector field is evaluated, we see that these vector fields must be tangent to $Y$. Because $\|s\|^2$ achieves a maximum in $Y$ at $y$, this means $\lied(\xi_{M}) \|s\|^2 (y) = 0$ for all $\xi \in \lieb^{\sigma}$, where $\xi_M$ is the fundamental vector field on $M$ induced by $\xi$.

Let $\func{\imaginary}{\lie^{\cc}}{\lie}$ be projection onto the imaginary component of $\lie^{\cc}$ with respect to the real form $\lie$, and let $\func{\pr}{\lieb}{\lie}$ be the restriction of $\imaginary$ to $\lieb \subset \lie^{\cc}$. For each $\xi \in \lie$ let $\func{\phi^{\xi}}{M}{\rr}$ be the function given by the pairing of $\Phi$ with elements of $\lie$: $\phi^{\xi} := \iprod{\xi}{\Phi}$. Then \cite[Equation 7]{guill} states that \[\lied(\xi_{M}) \|s\|^2 = 4 \pi r \left( - \lambda( \pr \xi) + \phi^{\pr \xi} \right) \|s\|^2\] for all $\xi \in \lieb$. By our reasoning in the previous paragraph, this tells us that \[0 = \lied(\xi_{M}) \|s\|^2 (y)= 4 \pi r \left( - \lambda( \pr \xi) + \phi^{\pr \xi}(y) \right) \|s(y)\|^2\] for all $\xi \in \lieb^{\sigma}$. Because $\|s(y)\|^2 > 0$, this implies that $- \lambda( \pr \xi) + \phi ( \pr \xi)(y) = 0$, and hence \begin{equation} \label{eq4} \iprod{\xi}{\Phi(y)} = \lambda (\pr \xi), \end{equation} for all $\xi \in \lieb^{\sigma}$. Recall that $\lambda \in \lieq^*$, and because $y \in M^{\tau}$ we also know $\Phi(y) \in \lieq^*$. Hence if we show that $\lieq \subset \pr(\lieb^{\sigma})$, then Equation~\ref{eq4} implies that $\Phi(y) = \lambda$.

Let $\varepsilon \in \lieq$. In \cite[page 10]{guill}, it is shown that $\func{\pr}{\lieb}{\lie}$ is onto. Therefore there exists some $\delta \in \lie$ such that $\delta + \ii \varepsilon \in \lieb$. Put $\zeta = \frac{1}{2} \left( \delta + \ii \varepsilon + \sigma(\delta+\ii \varepsilon) \right)$, and note that $\zeta$ is fixed by $\sigma$. Because $\lieb$ is preserved by $\sigma$ and is a vector space, we have $\zeta \in \lieb$. Since $\sigma$ is an extension of the involution $\func{\gamma}{\lie}{\lie}$, we know $\sigma(\delta) \in \lie$. One checks easily that $\imaginary(\zeta) = \varepsilon$, and therefore $\lieq \subset \pr(\lieb^{\sigma})$, and so $\Phi(y) = \lambda$. Thus $\mathcal{C}_{\gamma}(Y)$ is a subset of the rational points in $\Delta(Y)$.

\medskip

Now let $\lambda = \Phi(y) \in \Delta(Y)$ be a rational point. Since $Y \subset X$, $\Delta(Y) \subset \Delta(X)$, so $\lambda$ is a rational point of $\Delta(X)$ also. By Theorem~\ref{guillthm}, this means that $\lambda \in \mathcal{C}(X)$. Since $y \in Y \subset M^{\tau}$ we have $\lambda = \Phi(y) \in \lieq^*$. By Proposition~\ref{prop1}, $\lambda \in \mathcal{C}(X) \cap \lieq^* = \mathcal{C}_{\gamma}(Y)$. Thus the rational points of $\Delta(Y)$ are contained in $\mathcal{C}_{\gamma}(Y)$.
\end{proof}

We can now prove our main result, Theorem~\ref{mainthm}.

\begin{proof}[Proof of Theorem~\ref{mainthm}] We know that $\Delta(Y) \subset \Delta(X) \cap \lieq^*$. In the course of proving the main result of \cite{guill}, Guillemin and Sjamaar proved that $\Delta(X) = \overline{\mathcal{C}(X)}$, where the bar denotes the closure. Hence $\Delta(X) \cap \lieq^* = \overline{\mathcal{C}(X)} \cap \lieq^*$. Because $\lieq^*$ is equal to the closure of its rational points, as noted in Remark~\ref{rem1}, we know that $\overline{\mathcal{C}(X)} \cap \lieq^* = \overline{\mathcal{C}(X) \cap \lieq^*}$. Finally, Proposition~\ref{prop1} implies that $\overline{\mathcal{C}(X) \cap \lieq^*} = \overline{\mathcal{C}_{\gamma}(Y)}$. Therefore \begin{equation} \label{eq2} \Delta(X) \cap \lieq^* = \overline{\mathcal{C}_{\gamma}(Y)}. \end{equation}

Because $Y$ is a closed subset of the compact space $M$, it is compact. So $\Phi(Y)$ is compact in $\lie^*$ and hence closed. Therefore its intersection with the closed positive Weyl chamber, $\Delta(Y) = \Phi(Y) \cap \chamber$, is also closed. By Theorem~\ref{lem1} we know that $\mathcal{C}_{\gamma}(Y) \subset \Delta(Y)$, so $\overline{\mathcal{C}_{\gamma}(Y)} \subset \Delta(Y)$. Putting this together with Equation~\ref{eq2}, we see that $\Delta(X) \cap \lieq^* \subset \Delta(Y)$. Thus $\Delta(X) \cap \lieq^* = \Delta(Y)$.
\end{proof}

\section{Closures of Borel orbits}
\label{sec4}

Throughout this section, we will assume that the Borel subgroup $B$ is preserved by the antiholomorphic involution $\sigma$. 

The simplest example of a closed irreducible complex subvariety of $M$ preserved by $G^{\cc}$ is the closure of a $G^{\cc}$ orbit: $\overline{G^{\cc} m}$, for some $m \in M^{\tau}$. In \cite[Proposition 5.5]{oshea} it was shown that the ``real'' part of this subvariety, $(\overline{G^{\cc} m})^{\tau}$, has a nice decomposition. The simplest example of a closed irreducible complex subvariety of $M$ preserved by $B$ is the closure of a Borel orbit, and the ``real'' part of this subvariety has a corresponding decomposition. The proof is the same, after intersecting everything with $B$.

\begin{lemma}
Let $H$ denote the identity component of the real Lie group $B^{\sigma}$. For every $m \in M^{\tau}$, the set $(B m)^{\tau}$ has a finite number of components, each of which consists of a single $H$-orbit.
\end{lemma}

Therefore for any $m \in M$, $\overline{H m}$ is the closure of a component of $(\overline{B m})_{\text{reg}} \cap M^{\tau}$, so Theorem~\ref{mainthm} tells us that $\Delta ( (\overline{B m})^{\tau}) = \Delta(\overline{B m}) \cap \lieq^* = \Delta(\overline{H m})$.

Because our main result is so similar to that of Theorem~\ref{guillthm}, several of the corollaries of that theorem in \cite{guill} lead immediately to corresponding corollaries in our situation.

\begin{corollary} \label{cor1}
Suppose $B$ and $X$ are as in Theorem~\ref{mainthm}, and that $X^{\tau} \cap X_{\text{reg}} \neq \emptyset$. Then the set of $x \in X$ such that $\Delta(X^{\tau}) = \Delta( \overline{H x} )$ is nonempty and Zariski--open in $X$. Here $H$ is the identity component of $B^{\sigma}$.
\end{corollary}

\begin{proof}
\cite[Corollary 2.5]{guill} states that the set of $x \in X$ such that $\Delta(X) = \Delta(\overline{Bx})$ is nonempty and Zariski--dense in $X$. By Theorem~\ref{guillthm} this is equivalent to the statement that $\mathcal{C}(X) = \mathcal{C}(\overline{Bx})$, which in turn implies that $\mathcal{C}(X) \cap \lieq^* = \mathcal{C}(\overline{Bx}) \cap \lieq^*$. By Theorem~\ref{mainthm} and Corollary~\ref{maincor}, this means that $\Delta(X^{\tau}) = \Delta \left( (\overline{Bx})^{\tau}\right) = \Delta(\overline{Hx})$.
\end{proof}

\begin{corollary} \label{cor2}
The collection of polytopes $\Delta(X^{\tau})$, where $X$ ranges over all $B$ and $\tau$-invariant irreducible closed complex subvarieties of $M$, is finite.
\end{corollary}

\begin{proof}
In \cite[Corollary 2.6]{guill} it is proved that the collection of polytopes $\Delta(X)$, where $X$ ranges over the same set described in the statement of this corollary, is finite. Our corollary then follows immediately from the fact that each $\Delta(X^{\tau}) = \Delta(X) \cap \lieq^*$, by Theorem~\ref{mainthm}.
\end{proof}

Because $G^{\cc}$-invariance implies $B$-invariance, and because $M$ is itself both $B$ and $G^{\cc}$-invariant, Corollary~\ref{cor1} leads to the following result.

\begin{corollary}
Suppose $M^{\tau}$ contains a regular point. Then the set of $x \in M$ for which $\Delta( \overline{H x}) = \Delta( \overline{G' x}) = \Delta(M^{\tau})$ is nonempty and Zariski--open in $M$. Here $H$ is the identity component of $B^{\sigma}$ and $G'$ is the identity component of $(G^{\cc})^{\sigma}$.
\end{corollary}

\section{Examples}
\label{sec5}

Probably the most abundant source of examples to which the theorems in this paper apply is the constructions in the proof of the Borel--Weil Theorem. Suppose $G$ is a compact and connected Lie group, $T \subset G$ is a maximal torus, $\chamber \subset \liet^*$ is a choice of positive Weyl chamber, and $B$ is the Borel subgroup of $G^{\cc}$ corresponding to $\chamber$. Then for each choice of dominant weight $\lambda \in \chamber$, we can construct an integral, compact, and connected K\"ahler manifold in the form of a complex flag variety $M_{\lambda} := G^{\cc}/P_{\lambda}$,  and a holomorphic line bundle $L_{-\lambda}$ over $M$ whose Chern class is represented by $M_{\lambda}$'s K\"ahler form. (Here $P_{\lambda}$ is the parabolic subgroup of $G^{\cc}$ corresponding to $\lambda$.) Furthermore, the group $G$ acts on $(M_{\lambda},L_{-\lambda})$ in a natural K\"ahlerian fashion. A thorough treatment of this material can be found in Section 4.12 of \cite{duis}.

For any choice of Lie group involution $\gamma$ on $G$, so long as $\gamma$ preserves the Borel subgroup $B$ and the parabolic subgroup $P_{\lambda}$, we can easily construct involutions $\sigma$ on $G^{\cc}$, $\tau$ on $M_{\lambda}$, and $\beta$ on $L_{-\lambda}$ so that all of requirements described in Section~\ref{sec1} are satisfied. As in Section~\ref{sec4}, let $x \in M_{\lambda}^{\tau}$, let $X = \overline{B x}$, and let $Y$ be the closure of $H x$, where $H$ is the identity component of $B^{\sigma}$. So long as $M_{\lambda}^{\tau}$ is nonempty, we have a situation where we can apply all of the results of this paper. For an added twist, we can let $G = U \times U$ be the product of compact Lie groups. By pre--composing the action of $G$ on $M_{\lambda}$ with the diagonal map $U \to U \times U$, $u \mapsto (u,u)$, we obtain an action of $U$ on $M_{\lambda}$. It is well--known that this action is Hamiltonian with moment map obtained by post--composing the $G$-moment map with the projection $(\lieu \oplus \lieu)^* \cong \lieu^* \oplus \lieu^* \to \lieu^*$ defined by dualizing the diagonal map $\lieu \to \lieu \oplus \lieu$.

For a specific example, let $U = \mathbf{SU}(2)$, so that $U^{\cc} = \mathbf{SL}(2,\cc)$. Let $T$ consist of the diagonal matrices in $U$, and $B$ the upper triangular matrices in $U^{\cc}$. Then $N=[B,B]$ consists of the strictly upper triangular matrices in $U^{\cc}$. Define $\alpha \in \liet^*$ by $\begin{pmatrix} 2 \pi x \ii & 0 \\ 0 & - 2 \pi x \ii \end{pmatrix} \mapsto x \in \rr$. Then $\liet^* \cong \rr \cdot \alpha \subset \lie^*$, the positive Weyl chamber corresponding to $B$ is $\chamber = \rrnn \cdot \alpha$, and the weight lattice of $(U,T)$ is $\Lambda = \zz \cdot \alpha$, so $\Lambda_+ = \zznn \cdot \alpha$.

Now let $G = U \times U$, and take $T \times T$ as a maximal torus and $B \times B$ as a Borel subgroup of $G^{\cc} = U^{\cc} \times U^{\cc}$. The closed positive Weyl chamber of $\liet^* \oplus \liet^*$ is then $\liet_+^* \times \liet_+^*$, the weight lattice is $\Lambda \times \Lambda$, and the set of dominant weights is $\Lambda_+ \times \Lambda_+$. Let $\lambda_1, \lambda_2 \in \Lambda_+$ be nonzero dominant weights of $(U,T)$. Then $(\lambda_1, \lambda_2) \in \Lambda_+ \times \Lambda_+$ is a nonzero dominant weight of $(G, T \times T)$, and the corresponding parabolic subgroup is $B \times B$. Our flag variety in this case is $M:= G^{\cc}/(B \times B) = (U^{\cc} \times U^{\cc})/(B \times B) \cong (U^{\cc}/B) \times (U^{\cc}/B)$. Note that $U^{\cc}/B$ is isomorphic to $\comproj{1}$, two-dimensional complex projective space, so $M \cong \comproj{1} \times \comproj{1}$. Let $L_{(-\lambda_1,-\lambda_2)}$ denote the holomorphic line bundle over $M$ constructed from the Borel--Weil Theorem. For involutions on $U$, $M$, and $L$, we can take standard complex conjugation. For $U=G \times G$ this means conjugation on each factor. It is easily verified that these satisfy all of the necessary compatibility conditions.

Let $r \in \zplus$. Then $(r \lambda_1, r \lambda_2)$ is also a dominant weight, so we can repeat the above construction, but we simply have $L_{(-r \lambda_1, -r \lambda_2)} \approx \left( L_{(-\lambda_1,-\lambda_2)} \right)^r$. Note also that \begin{equation} \label{eq4} L \approx L_{-\lambda_1} \boxtimes L_{-\lambda_2} \text{ and } L^r \approx L_{-r\lambda_1} \boxtimes L_{-r\lambda_2}. \end{equation} Recall that the space $\Gamma(M,L^r)$ of global holomorphic sections of $(M,L^r)$ is isomorphic as a $G$-representation to $V(r \lambda_1,r \lambda_2)^*$, the dual of the irreducible representation of $G$ with highest weight $(r \lambda_1,r \lambda_2)$. Similarly, $\Gamma(\comproj{1}, L_{-r\lambda_1}) \cong V(r\lambda_1)^*$ and $\Gamma(\comproj{1}, L_{-r\lambda_2}) \cong V(r\lambda_2)^*$. In general, we have the formulas $V(r \lambda_1)^* \cong V(-w_0 r \lambda_1)$ and $V(r \lambda_2)^* \cong V(-w_0 r \lambda_2)$, where $w_0$ is the longest element of the Weyl group of $(U,T)$. For $U=\mathbf{SU}(2)$, $w_0$ is the identity. Together with the equalities in (\ref{eq4}), this implies that $\Gamma(M,L^r) = \Gamma (\comproj{1}, L_{-r\lambda_1}) \otimes \Gamma (\comproj{1}, L_{-r\lambda_1}) \cong V(r \lambda_1) \otimes V(r \lambda_2)$.

If $\lambda \in \Lambda_+$ is a dominant weight of $(U,T)$, then $V_{\lambda}$ is equivalent to the space of homogeneous complex polynomials of degree $\lambda$ in two variables. (See pp. 305-306 of \cite{duis}.) If $F$ is such a polynomial and $u \in U$, then $(u \cdot F)(x,y) = F \left( u^{-1} (x,y) \right)$, where $U$ acts on $\ccc{2}$ in the usual way. Using this description, we see that $\Gamma(M,L^r) \cong V(r \lambda_1) \otimes V(r \lambda_2)$ can be viewed as the space of complex polynomials $F(x_1,y_1,x_2,y_2)$ which are homogenous of degree $r \lambda_1$ in the first two variables and homogeneous of degree $r \lambda_2$ in the last two variables, which is a vector space of dimension $(r \lambda + 1)(r \lambda_2 + 1)$.

By the Clebsch-Gordan formula, the dominant weights that appear as highest weights in the decomposition of $V(r \lambda_1) \otimes V(r \lambda_2)$ into irreducible representations are exactly $r(\lambda_1 + \lambda_2) - 2k$ for integers $k = 0, \ldots, \min \{r\lambda_1, r\lambda_2\}$. For each of these weights, there is a one-dimensional subspace of $V(r \lambda_1) \otimes V(r \lambda_2)$ which is $N$-invariant and on which $T$ acts by the given weight. Some careful computation shows that, for each $k = 0, \ldots, \min\{r\lambda_1,r\lambda_2\}$, this one-dimensional subspace is the complex span of the polynomial \begin{eqnarray*} F_{r,k} (x_1,y_1,x_2,y_2) &:=& \sum_{j=0}^k (-1)^{k-j} {k \choose j} x_1^j y_1^{r \lambda_1-j} x_2^{k-j} y_2^{r\lambda_2-k+j} \\ &=& y_1^{r\lambda_1-k} y_2^{r \lambda_2 - k} \left( x_1 y_2 - x_2 y_1 \right)^k.\end{eqnarray*} (The last equality follows from the Binomial Expansion Theorem.) From the first formula, it is easy to see that multiples of $F_{r,k}$ transform under $T$ according to the weight $r(\lambda_1+\lambda_2)-2k$, while the second formula allows easy verification that $F_{r,k}$ and its multiples are invariant under $N$. 

Fix $x=\left( (a_1:c_1),(a_2:c_2) \right) \in \comproj{1} \times \comproj{1}$, and put $X = \overline{B x}$. Notice that $F_{r,k} (b x)$ is a nonzero multiple of $F_{r,k} (x)$ for any $b \in B$, so $F_{r,k}$ vanishes on $X$ exactly when $F_{r,k}(x)=0$. Depending on the specific value of $x$, $X$ can be all of $M$, $\comproj{1}$ embedded diagonally in $M$, the product of $\comproj{1}$ and the point $(1:0)$ (in either order), or just $\left((1:0),(1:0)\right)$. Using our work above and considering the different possibilities zero versus nonzero possibilities for the homogeneous coordinates of $x$, one can calculate the highest weight polytope $\mathcal{C}(X)$, and then determine moment polytope $\Delta(X) \subset \liet^*$. For $X=M$, the polytope is the closed line segment $\left[ |\lambda_1-\lambda_2|, \lambda_1+\lambda_2 \right]$. For the diagonal of $M$, the polytope is the point $\{\lambda_1+\lambda_2\}$. For $X = \comproj{1} \times \{(1:0)\}$ or $X = \{(1:0)\} \times \comproj{1}$, one of the polytopes will be empty and the other will be the point $\{|\lambda_1 - \lambda_2|\}$, the order depending on the value of $\min \{ \lambda_1, \lambda_2\}$. For $X = \{ \left( (1:0), (1:0) \right) \}$, the polytope is empty.

Suppose $X^{\tau}$ is nonempty. If $X$ is $M$, the diagonal in $M$, the product of $\comproj{1}$ and a point, or just a point, then $X^{\tau}$ is $\reproj{1} \times \reproj{1}$, $\reproj{1}$ embedded in $\comproj{1}$ diagonally, the product of $\reproj{1}$ and a point, or just a point, respectively. Let $Y = X^{\tau}$. Because $\lie = \mathfrak{su}(2)$, all entries of elements of $\liet$ are pure imaginary, so $\gamma$ acts on $\liet^*$ by negation, so $\liet^* \subset \lieq^*$. Then Theorem~\ref{mainthm} implies that $\Delta(Y) = \Delta(X) \cap \lieq^* = \Delta(X)$, so that $\Delta(Y)$ is always the same polytope as $\Delta(X)$.

One can calculate that the closure of any $G^{\cc}$-orbit in $M$ is either all of $M$ or $\comproj{1}$ embedded diagonally in $M$, the real parts of which are $\reproj{1} \times \reproj{1}$ and $\reproj{1}$ embedded diagonally in $M$, respectively. The corresponding moment polytopes are the closed line segment $\left[ |\lambda_1-\lambda_2|, \lambda_1+\lambda_2 \right]$ and the point $\{ \lambda_1 + \lambda_2 \}$, respectively. Note that three of the closures of Borel orbits described above are not preserved by $G^{\cc}$, and for the Borel case we obtain an additional possible moment polytope (or two additional ones, if we include the empty set). So even in this relatively easy example, we find some situations to which the results of O'Shea and Sjamaar do not apply, while the new results of this paper do.

\section*{Acknowledgements}
The author would like to thank his family and friends for their unwavering support, and Professor Reyer Sjamaar for his excellent and crucial advice.

\end{document}